\documentclass[authoryear,8pt]{elsarticle}
\usepackage[dvips]{epsfig}
\usepackage{latexsym}
\usepackage{amsmath,amsthm}
\usepackage{natbib}
\usepackage{amssymb}
\usepackage{graphicx}
\usepackage[latin1]{inputenc}
\usepackage{bm}
\usepackage[notablist,nofiglist,tabhead]{endfloat}

\usepackage{geometry}
\def\R{\mathbb{R}}

\def\d{{\rm d}}
\def\1{\mathbf{1}}

\def\cov{{\rm cov}}

\def\C{\mathbb{C}}
\def\B{\mathbb{B}}
\def\D{\mathbb{D}}

\def\Z{\mathbb{Z}}

\newcommand{\ip}[1]{\lfloor #1 \rfloor}
\renewcommand{\vec}{\bm}

\newcommand{\p}{\overset{\Pr}{\to}}

\newcommand{\bstar}{m}

\newtheorem{prop}{Proposition}[section]

\newtheorem{cor}{Corollary}[section]
\newtheorem{lemma}{Lemma}[section]
\newtheorem{definition}{Definition}[section]
\newtheorem{condition}{Condition}[section]

\newdefinition{rmk}{Remark}

\title{Some results on change-point detection in cross-sectional dependence of multivariate data with changes in marginal distributions.}
\author{Tom Rohmer$^{1,2}$\\ \small{$^1$Universit\'e de Bordeaux, ISPED, Bordeaux France} \\
\small{$^2$Inserm,  Centre Inserm U-1219, Bordeaux France}\\
 146 rue L\'eo Saignat, 33000 Bordeaux.\\ \small{tom.rohmer@isped.u-bordeaux2.fr}}
\begin{document}


\begin{abstract}
Tests for break points detection in the law of random vectors have been proposed in several papers. Nevertheless, they have often little power for alternatives involving a change in the dependence between components of vectors. Specific tests for detection of a change in the copula of random vectors have also been proposed in recent papers, but they do not allow to conclude of a change in the dependence structure without condition that the margins are constant. The goal of this article is to propose a test for detection of a break in the copula when changes in marginal distribution occurs at  known instants. The performances of this test are illustrated by Monte Carlo simulations.\\
\end{abstract}
\begin{keyword}
non-parametric tests, sequential empirical copula process, Monte Carlo experiments
\end{keyword}

\maketitle

\section{Introduction}

Let $\vec X$ be a $d$-dimensional random vector ($d\geq 2$), with cumulative distribution function (c.d.f.) $F$ and marginal cumulative distribution functions (m.c.d.f.s) $F_1,\ldots,F_d$.
When the m.c.d.f.s $F_1,\ldots,F_d$ are continuous, Sklar's Theorem \citep[see][]{Skl59} allows us to say that there exists a unique function $C$ called copula, characterizing the dependence of random vector $\vec X$, such that $F$ can be written as:
\begin{equation}\label{Sklar}
F(\vec x) = C(F_1(x_1),\ldots,F_d(X_d)),\qquad \vec x\in\R^d.
\end{equation}
Let $\vec X_1,\ldots,\vec X_n$ be $d$-dimensional observations.
The purpose of change-points detection is to test the hypothesis 
\begin{equation}\label{H'0}
H_0:\exists F \text{ such as }\vec X_1,\ldots,\vec X_n \text{  have c.d.f. } F,
\end{equation}
against $\neg H_0$. The equation~(\ref{Sklar}) involves that $H_0$ can be rewritten as $H_0 = H_{0,m}\cap H_{0,c}$, with
\begin{align}
\label{H0m}H_{0,m}&: \exists F_1,\ldots F_d\text{ such as }  \vec X_1,\ldots,\vec X_n \text{  have m.c.d.f. } F_1,\ldots F_d,\\
\label{H0c}H_{0,c}&: \exists C\text{ such as }  \vec X_1,\ldots,\vec X_n \text{  have copula } C.
\end{align}
A change either in the copula of random vectors or in one of m.c.d.f.s implies the rejection of the null hypothesis $H_0$. Many non-parametric tests for $H_0$ based on empirical processes are present in the literature; see for example \citet{Bai94,CsoHor97,Ino01}. 
These tests are not very sensitive to detect a change in  the copula which leaves the m.c.d.f.s unchanged. This conclusion is highlighted in \citet[Section 2]{HolKojQue13} through Monte Carlo simulations.\\


Non-parametric tests for break detection, sensitive to changes in the copula of observations and based on the two-sided sequential empirical copula process are considered in \cite{BucKojRohSeg14}. 
These tests do not allow to conclude in favour of $\neg H_{0,c}$ without condition on the constancy of m.c.d.f.s.
In many situations, see for example section~5.2, a specific event can lead to changes in the marginal cumulative distributions. The question then becomes whether the specific event changes the copula or not. 
The aim of this paper is to propose a test to detect a change in the dependence structure of random vectors, sensitive to changes in copula of observations and adapted in the case of alternative hypotheses involving abrupt changes in m.c.d.f.s.\\  

This paper is organized as follows. The procedure to test the null hypothesis of a break in c.d.f.\ when a change in the m.c.d.f.s occurs is presented in Section~2. An adaptation of results of Section~2 when multiple changes in m.c.d.f.s occur is described in Section~3. Section~4 contains the results of Monte Carlo simulations. Finally, Section~5 reports brief discussions about the case of $alpha$-mixing observations and presents an illustration on a specific situation.

\section{Break detection in the copula when a break time in the m.c.d.f.s is known}

In the sequel, the weak convergence, denoted by $\leadsto$, must be understood as being the weak convergence in the sense of Definition 1.3.3 in \cite{vanWel96}. For a set $T$, $\ell^{\infty}(T)$ denotes the space of bounded real-valued functions on $T$ equipped with the uniform metric. \\

 Let $\vec X_1,\ldots,\vec X_n$ be $d$-dimensional random vectors ($d\geq 2$) and consider for $1\leq k\leq l\leq n$ the empirical copula $C_{k:l}$ of the sub-sample $\vec X_k,\ldots,\vec X_l$ as suggested in \cite{Deh79}:
\begin{equation}\label{cp}
C_{k:l}(\vec u)=\frac{1}{l-k+1}\sum_{i=k}^{l}\prod_{j=1}^d\1(F_{k:l,j}(X_{ij})\le u_j),\qquad \vec u\in[0,1]^d,
\end{equation}
where for $j=1,\ldots,d$, $F_{k:l,j}$ is the empirical cumulative distribution function (e.c.d.f.) of sample $X_{kj},\ldots,X_{lj}$:
\begin{equation}\label{fdre}
F_{k:l,j}(x) = \frac{1}{l-k+1}\sum_{i=k}^{l}\1(X_{ij}\le x),\qquad x\in\R.
\end{equation}
In \cite{BucKojRohSeg14}, the following Cram\'er--von Mises's type statistic to test $\neg H_0$ is suggested:
\begin{equation}\label{Sn}
S_n = \sup_{s\in[0,1]}\sqrt{n}\lambda_n(s,1)\lambda_n(0,s)\int_{[0,1]^d}\{C_{1:\ip{ns}}(\vec u) - C_{\ip{ns}+1:n}(\vec u)\}^2\d C_{1:n}(\vec u),
\end{equation}
where $\lambda_n(s,t)=(\ip{nt}-\ip{ns})/n$, $s\leq t\in[0,1]$. \\

Monte Carlo simulations \citep[see section 5 of][]{BucKojRohSeg14} highlighted that a strategy of boostraping with independent or dependent multipliers according to the observations  \citep[see ][]{BucKoj14, BucKojRohSeg14} of the statistic $S_n$ leads to very good performances in term of powers for alternatives hypotheses that involve a change in copula which leave the m.c.d.f.s\ unchanged.\\

Let us suppose that it exists a break time $\bstar=\ip{nb}$ in m.c.d.f.s, $b\in(0,1)$ known. We propose a test for $H^m_0 = H_{0,c}\cap H_{1,m}$,  where  $H_{0,c}$ is defined in \eqref{H0c} and $H_{1,m}$ is defined by:
\begin{equation}\label{H1mb}
H_{1,m}: \exists F_1,\ldots F_d \text{ and } F'_1,\ldots F'_d \text{ such that } \begin{matrix}\vec X_1,\ldots,\vec X_{\bstar}\text{  have m.c.d.f. } F_1,\ldots F_d,\\
\vec X_{\bstar+1},\ldots,\vec X_n\text{  have m.c.d.f. } F'_1,\ldots F'_d. \end{matrix}
\end{equation}
Note that we do not suppose that $F'_1,\ldots F'_d$ are necessarily different from $F_1,\ldots F_d$.
In other words, we do not assume a change in the m.c.d.f.s. However we suppose that if there is a change in the m.c.d.f.s, it is a unique and abrupt change at time $\bstar$. 

Let $\vec X_1,\ldots,\vec X_n$ be $d$-dimensional random vectors with unknown copula $C$, such that $\vec X_1,\ldots,\vec X_{\bstar}$ have m.c.d.f.s $F_1,\ldots,F_d$ and $\vec X_{\bstar+1},\ldots,\vec X_n$ have m.c.d.f.s $F'_1,\ldots,F'_d$, where $F_1,\ldots,F_d, F'_1,\ldots,F'_d$ are unknown and the break time $\bstar=\ip{nb}$, $b\in(0,1)$ is known.\\

For $i\in\{1,\ldots,n\}$, let us consider the random vectors $\vec U_{i,\bstar}$ defined by
\begin{equation}\label{Uim}
\vec U_{i,\bstar}=\left\{\begin{aligned}[lr](F_1(X_{i1}),\ldots,F_d(X_{id})) &\hspace{0.5cm} i\in\{1,\ldots,\bstar\}\\(F'_1(X_{i1}),\ldots,F'_d(X_{id})) &\hspace{0.5cm} i\in\{\bstar+1,\ldots,n\}. \end{aligned}\right. 
\end{equation}

Note that the vectors $\vec U_{i,\bstar}$, $i=1,\ldots,n$ have $C$ for c.d.f. For $i\in\{1,\ldots,n\}$, let $\hat{\vec U}^{1:n}_{i,\bstar}$ defined by:

$$
\hat{\vec U}^{1:n}_{i,\bstar}=\left\{\begin{aligned}[lr](F_{1:\bstar,1}(X_{i1}),\ldots,F_{1:\bstar,d}(X_{id})) &\hspace{0.5cm} i\in\{1,\ldots,\bstar\}\\(F_{\bstar+1:n,1}(X_{i1}),\ldots,F_{\bstar+1:n,d}(X_{id})) &\hspace{0.5cm} i\in\{\bstar+1,\ldots,n\}, \end{aligned}\right. 
$$
where for all $1\leq k\leq l \leq n$ and $j=1,\ldots,d$ $F_{k:l,j}$ is the empirical c.d.f. of $X_{kj},\ldots,X_{lj}$ as defined in~(\ref{fdre}).
The vectors $\hat{\vec U}^{1:n}_{i,\bstar}$, $i=1,\ldots,n$ can be seen as pseudo-observations of copula $C$. 
An estimator of $C$ is given by the empirical distribution of $\hat{\vec U}^{1:n}_{1,\bstar},\ldots,\hat{\vec U}^{1:n}_{n,\bstar}$:
$$
C_{1:n,\bstar}(\vec u) = \frac{1}{n}\sum_{i=1}^n\1(\hat{\vec U}^{1:n}_{i,\bstar}\le \vec u),\quad \vec u\in[0,1]^d.
$$
This estimator can be rewritten for $\vec u\in[0,1]^d$ by:
\begin{align*}
C_{1:n,\bstar}(\vec u) 
&= \frac{\bstar}{n}C_{1:\bstar}(\vec u) + \frac{n-\bstar}{n}C_{\bstar+1:n}(\vec u),
\end{align*}
where for any subsample $\vec X_k,\ldots,\vec X_l$, $1\leq k\leq l\leq n$, $C_{k:l}$ is the empirical c.d.f.  of random vectors $\hat{\vec U}^{k:l}_k,\ldots,\hat{\vec U}^{k:l}_l$ defined in~(\ref{cp}), and for $i=k,\ldots,l$
$$
\hat{\vec U}^{k:l}_i = (F_{k:l,1}(X_{i1}),\ldots,F_{k:l,d}(X_{i1})).
$$
More particularly, for $1\le k\le l\le n$, $\bstar\in\{1,\ldots,n-1\}$ and $\vec u\in[0,1]^d$,
\begin{equation}\label{Cklm}
C_{k:l,\bstar}(\vec u) = \left\{\begin{array}{ll}\frac{\bstar-k+1}{l-k+1}C_{k:\bstar}(\vec u) + \frac{l-\bstar}{l-k+1}C_{\bstar+1:l}(\vec u) &\quad \bstar\in[k,l],\\C_{k:l}(\vec u) &\quad \bstar\notin[k,l].\end{array}\right. 
\end{equation}

For a subsample $\vec X_k\ldots,\vec X_l$, $1\leq k\leq l\leq n$ consider the following pseudo-observations of copula $C$:
\begin{equation*}
\hat{\vec U}^{k:l}_{i,\bstar}=\left\{
\begin{aligned}
&\left\{
\begin{aligned}[lr](F_{k:\bstar,1}(X_{i1}),\ldots,F_{k:\bstar,d}(X_{i1})) &\hspace{0.5cm} i\in\{k,\ldots,\bstar\}\\(F_{\bstar+1:l,1}(X_{i1}),\ldots,F_{\bstar+1:l,d}(X_{i1})) &\hspace{0.5cm} i\in\{\bstar+1,\ldots,l\} \end{aligned}
\right. &\hspace{0.5cm} m\in[k,l],\\
&(F_{k:l,1}(X_{i1}),\ldots,F_{k:l,d}(X_{i1}))&\hspace{0.5cm} m\notin[k,l],
\end{aligned}
\right.
\end{equation*}
$F_{k:l}$ defined in \eqref{fdre}. Then $C_{k:l,m}$ defined in \eqref{Cklm} is the empirical cumulative distribution of $\hat{\vec U}^{k:l}_{k,\bstar},\ldots,\hat{\vec U}^{k:l}_{l,\bstar}$. \\

The corresponding two-sided sequential empirical copula process is defined by

\begin{align}\label{Cnm}
\C_{n,m}(s,t,\vec u) &= \sqrt{n}\lambda_n(s,t)\{C_{\ip{ns}+1:\ip{nt},m}(\vec u) - C(\vec u)\},\quad (s,t,\vec u)\in\Delta\times[0,1]^d,\\
\nonumber&=\frac{1}{\sqrt{n}}\sum_{i=\ip{ns}+1}^{\ip{nt}}\{\1(\hat{\vec U}^{\ip{ns}+1:\ip{nt}}_{i,m}\leq \vec u) - C(\vec u)\},
\end{align}
$\Delta=\{(s,t)\in[0,1]^2|s\leq t\}$. The test statistic proposed in this paper is based on the process $\D_{n,\bstar}$, defined by
$$
\D_{n,\bstar}(s,\vec u) = \sqrt{n}\lambda_n(0,s)\lambda_n(s,1)\{C_{1:\ip{ns},\bstar}(\vec u) - C_{\ip{ns}+1:n,\bstar}(\vec u)\},\quad (s,\vec u)\in[0,1]^{d+1}.
$$

 Note that $\D_{n,m}$ can be rewritten as
$$
\D_{n,m}(s,\vec u) = \lambda_n(s,1)\C_{n,m}(0,s,\vec u) - \lambda_n(0,s)\C_{n,m}(s,1,\vec u),\qquad (s,\vec u)\in[0,1]^{d+1}.
$$

Similarly to $S_n$ defined in~\eqref{Sn}, we consider the Cram\'er--von Mises statistic
\begin{align}\label{Snb}
S_{n,\bstar} &= \sup_{s\in[0,1]}\int_{[0,1]^d}\{\D_{n,\bstar}(s,\vec u)\}^2\d C_{1:n,\bstar}(\vec u)\\
\nonumber&=\max_{1\leq k\leq n-1}\frac{1}{n}\sum_{i=1}^n\{\D_{n,\bstar}(k/n,\hat{\vec U}^{1:n}_{i,\bstar})\}^2.
\end{align}

The asymptotic behaviour of the empirical process $\D_{n,\bstar}$ is given on Proposition~\ref{prop1}, proved in Appendix A. The result is obtained under the following non-restrictive condition, proposed in \citet{Seg12}:

\begin{condition}\label{Cj1}
For any $j\in \{1,\ldots,d\}$, the partial derivatives $\dot C_j = \partial C/\partial u_j $ exist and are continuous on $V_j=\{\vec u\in[0,1]^d,u_j\in(0,1)\}$.
\end{condition}

\begin{prop}\label{prop1}
Let $\vec X_1,\ldots,\vec X_n$ be $d$-dimensional independent random vectors with copula $C$, such that for $b\in(0,1)$ known and $m=\ip{nb}$, the random vectors $\vec X_1,\ldots,\vec X_{\bstar}$ have m.c.d.f.s $F_1,\ldots,F_d$ and the random vectors $\vec X_{\bstar+1},\ldots,\vec X_n$ have m.c.d.f.s $F'_1,\ldots,F'_d$. \\


Then, under Condition \ref{Cj1}, the process $\D_{n,\bstar}$ converges weakly in $\ell^{\infty}([0,1]^{d+1})$, to a stochastic process $\D_C$ defined by
\begin{equation}\label{DC}
\D_C(s,\vec u) = \C^0_C(s,\vec u) - s\C^0_C(1,\vec u),\quad (s,\vec u)\in[0,1]^{d+1},
\end{equation}
where for $(s,\vec u)\in[0,1]^{d+1}$,
\begin{equation}\label{Cc}
\C^0_C(s,\vec u) = \Z_C(s,\vec u)-\sum_{j=1}^d\dot C_j(\vec u)\Z_C(s,\vec u^{\{j\}})
\end{equation}
with $\Z_C$ is a tight centred Gaussian process with covariance function
$$
\cov\{\Z_C(s,\vec u),\Z_C(t,\vec v)\} = \min(s,t)\{C(\vec u\wedge \vec v) - C(\vec u)C(\vec v)\},\quad (s,\vec u),(t,\vec v)\in[0,1]^{d+1} ,
$$
$\vec u\wedge \vec v = (\min(u_1,v_1),\ldots,\min(u_d,v_d))$ and  $\vec u^{\{j\}} = (1,\ldots,1,u_j,1,\ldots,1)$.
\end{prop}
To resample $\D_{n,m}$, we note that for $b\leq s$ and $\vec u\in[0,1]^d$, $\D_{n,m}$ rewrites as
$$
\D_{n,m}(s,\vec u) =  \lambda_n(s,1)\{\C_{n,m}(0,b,\vec u) + \C_{n,m}(b,s,\vec u)\} - \lambda_n(0,s)\C_{n,m}(s,1,\vec u)
$$
and for $b\geq s$ and $\vec u\in[0,1]^d$ as
$$
\D_{n,m}(s,\vec u) =  \lambda_n(s,1)\C_{n,m}(0,s,\vec u) - \lambda_n(0,s)\{\C_{n,m}(s,b,\vec u) + \C_{n,m}(b,1,\vec u)\}.
$$

Let B a large integer and consider for $(s,t,\vec u)\in \Delta\times[0,1]^{d}$ and for $\beta=1,\ldots,B$, the processes
$$
\check\B^{(\beta)}_{n,m}(s,t,\vec u)=\frac{1}{\sqrt{n}}\sum_{i=\ip{ns}+1}^{\ip{nt}}\xi^{(\beta)}_{i}\left\{\1(\hat{\vec U}^{\ip{ns}+1:\ip{nt}}_{i,m}\leq \vec u) - C_{\ip{ns}+1:\ip{nt},m}(\vec u)\right\},
$$
and 
\begin{equation}\label{checkCn}
\check{\C}^{(\beta)}_{n,m}(s,t,\vec u) = \check\B^{(\beta)}_{n,m}(s,t,\vec u)- \frac{1}{\sqrt{n}}\sum_{j=1}^d \dot C_{j,\ip{ns}+1:\ip{nt},m}(\vec u)\check\B^{(\beta)}_{n,m}(s,t,\vec u^{\{j\}}),
\end{equation}
with for $\beta=1, \ldots,B$, and $1\leq i\leq n$, $\xi^{(\beta)}_{i}$ are i.i.d. standard normal random variables.\\

Re-sampling versions $\check{\D}_{n,\bstar}^{(\beta)}$ of  $\D_{n,\bstar}$ can be obtained for $\beta=1,\ldots,B$ and $(s,\vec u)\in [b,1]\times[0,1]^d$ by
$$
\check{\D}^{(\beta)}_{n,m}(s,\vec u) =  \lambda_n(s,1)\{\check{\C}^{(\beta)}_{n,m}(0,b,\vec u) + \check{\C}^{(\beta)}_{n,m}(b,s,\vec u)\} - \lambda_n(0,s)\check{\C}^{(\beta)}_{n,m}(s,1,\vec u)
$$
and for $(s,\vec u)\in [0,b]\times[0,1]^d$ by
$$
\check{\D}^{(\beta)}_{n,m}(s,\vec u) =  \lambda_n(s,1)\check{\C}^{(\beta)}_{n,m}(0,s,\vec u) - \lambda_n(0,s)\{\check{\C}^{(\beta)}_{n,m}(s,b,\vec u) + \check{\C}^{(\beta)}_{n,m}(b,1,\vec u)\}.
$$

For $j=1,\ldots,d$, the functions $\dot C_{j,\ip{ns}+1:\ip{nt},m}$ appearing in \eqref{checkCn} are an adaptation of the estimator of $\dot C_j$ proposed in section 4.2 in~\cite{BucKojRohSeg14} consisting in simple differencing at a bandwidth $h_{k:l} = \min\{(l-k+1)^{-1/2},1/2\}$ of the empirical copula process:

$$
\dot C_{j,k:l,m}(\vec u) = \frac{C_{k:l,m}(\vec u^{j,+}) - C_{k:l,m}(\vec u^{j,-})}{u_j^+ - u_j^-},\qquad \vec u \in[0,1]^d,\quad 1\leq k\leq l\leq n,
$$
with $u_j^+ = \min(u_j + h_{k,l},1)$, $u_j^- = \max(u_j - h_{k,l},0)$ and $\vec u^{j,\pm} = (u_1,\ldots,u_j^{\pm},\ldots,u_d)$. This estimator is in spirit of section~3 of \cite{KojSegYan11}.\\

Note that $\check{\D}^{(\beta)}_{n,m}$ is slightly different to 
$$
\widetilde{\D}^{(\beta)}_{n,m}(s,\vec u) = \lambda_n(s,1)\check{\C}^{(\beta)}_{n,m}(0,s,\vec u) - \lambda_n(0,s)\check{\C}^{(\beta)}_{n,m}(s,1,\vec u),\qquad (s,\vec u)\in[0,1]^{d+1}.
$$
These resample $\widetilde{\D}^{(\beta)}_{n,m}$ can be studied in a future research.\\

We have the following Proposition (proved in Appendix A.) 

\begin{prop}\label{prop:Cnmboot}
Under the same condition as Proposition~\ref{prop1}, we have the following result:
$$
\left(\D_{n,m},\check\D^{(1)}_{n,m},\ldots,\check\D^{(B)}_{n,m}\right)\leadsto \left(\D_C, \D^{(1)}_C,\ldots,\D^{(B)}_C\right), 
$$ 
in $\ell^{\infty}(\Delta\times[0,1]^d)^{B+1}$, where for $(s,t,\vec u)\in\Delta\times[0,1]^d$, where $\D_C$ is defined in \eqref{DC}  and $\D^{(1)}_C,\ldots,\D^{(B)}_C$ are independent copies of $\C_C$.
\end{prop}

As a corollary of Proposition~\ref{prop2} and continuous mapping theorem, we have the following result:
\begin{cor}
Let $\vec X_1,\ldots,\vec X_n$ be $d$-dimensional independent random vectors with copula $C$, such that for $b\in(0,1)$ known and $m=\ip{nb}$, the random vectors $\vec X_1,\ldots,\vec X_{\bstar}$ have c.d.f. $F$ and the random vectors $\vec X_{\bstar+1},\ldots,\vec X_n$ have c.d.f. $F'$. \\

Consider the statistic defined in \eqref{Snb} by
$$
S_{n,\bstar} = \sup_{s\in[0,1]}\int_{[0,1]^d}\{\D_{n,\bstar}(s,\vec u)\}^2\d C_{1:n,\bstar}(\vec u),
$$
and re-sampling versions of this statistic defined for $\beta=1,\ldots,B$ by
$$
\check S^{(\beta)}_{n,\bstar} = \sup_{s\in[0,1]}\int_{[0,1]^d}\{\check{\D}^{(\beta)}_{n,\bstar}(s,\vec u)\}^2\d C_{1:n,\bstar}(\vec u).
$$
Under the Condition \ref{Cj1},
$$
(S_{n,\bstar},\check S^{(1)}_{n,\bstar},\ldots,\check S^{(B)}_{n,\bstar})\leadsto (S_C,S^{(1)}_C,\ldots,S^{(B)}_C),
$$
where
$$
S_C = \sup_{s\in[0,1]}\int_{[0,1]^d}\{\D_C(s,\vec u)\}^2\d C(\vec u),
$$
and $S^{(1)}_C,\ldots,S^{(B)}_C$ are independent copies of $S_C$.
\end{cor}

An approximate p-value of the test for $H^m_0$ can be obtained by
$$
\hat{p}^m_{n,B} = \frac{1}{B}\sum_{\beta=1}^B\1\left\{S^{(\beta)}_{n,\bstar} \geq S_{n,\bstar}\right\}.
$$ 
The previous proposition and the Proposition~F.1 in the supplementary material of \cite{BucKoj14} allow to conclude that the test based on $\hat{p}_{n,\bstar,B}$ will hold its level asymptotically as $n\rightarrow \infty$ followed by $B \rightarrow \infty$.

\section{Break detection in the copula when multiple break times in m.c.d.f.s are known }

Suppose in this section that for an integer $p>0$ and for $j=1,\ldots,p+1$, the random vectors $\vec X_{m_{j-1}+1},\ldots, \vec X_{m_{j}}$ have m.c.d.f.s $F_{1j},\ldots,F_{dj}$ where $m_0 = 0$, $m_{p+1} = n$ and for $j=1,\ldots,p$ $m_j=\ip{nb_j}$, $0\leq b_1 < \ldots < b_p\leq 1$ known. In the sequel, $\vec m$ denotes the vector of break points $(m_1,\ldots,m_p)$. Similarly at the section~2, we propose a test for $H_0^{\vec m} = H_{0,c}\cap H_{1,\vec m}$ where $H_{1,\vec m}$ is defined by
$$
H_{1,\vec m} : \text{for }j=1,\ldots,p\text{ }\exists F_{1j},\ldots,F_{dj}\text{ such that }  \vec X_{m_{j-1}+1},\ldots,\vec X_{m_{j}} \text{ have m.c.d.f.s }F_{1j},\ldots,F_{dj}.
$$

Here, for $i=1,\ldots,n$ we consider the random vectors $\hat{\vec U}^{1:n}_{i,\vec m}$ defined by
$$
\hat{\vec U}^{1:n}_{i,\vec m} = \left\{
\begin{matrix}
\left(F_{1:m_1,1}(X_{i1}),\ldots,F_{1:m_1,d}(X_{id})\right) &  m_0+1=1\leq i \leq m_1\\
\left(F_{m_1+1:m_2,1}(X_{i1}),\ldots,F_{m_1+1:m_2,d}(X_{id})\right) &  m_1< i \leq m_2\\
\vdots & \vdots\\
\left(F_{m_p+1:n,1}(X_{i1}),\ldots,F_{m_p+1:n,d}(X_{id})\right) &  m_p< i \leq m_{p+1} = n,
\end{matrix}\right.
$$
where for $j=1,\ldots,p+1$ and $q=1,\ldots,d$, $F_{m_{j-1}+1:m_j,q}$ is the empirical c.d.f. of $X_{m_{j-1}+1q},\ldots,X_{m_{j}q}$ defined in~\eqref{fdre}.
For $1\leq k\leq l\leq n$, denote by $C_{k:l,\vec m}$ the empirical c.d.f. of random vectors $\hat{\vec U}_{1,\vec m},\ldots,\hat{\vec U}_{n,\vec m}$.
More particularly, for $1\le k\le l\le n$, $1<m_1<\ldots<m_p<n$ and $\vec u\in[0,1]^d$, we have
$$
C_{k:l,\vec m}(\vec u) = \left\{\begin{array}{ll}C_{k:l}(\vec u) &\quad m_1,\ldots,m_p\notin[k,l] \\
\dfrac{1}{l-k+1}\displaystyle\sum_{j=q_1}^{q_{2+1}}(m'_j - m'_{j-1})C_{m'_{j-1}+1:m'_j}(\vec u) & \quad \underset{(m'_{q_1},\ldots,m'_{q_2})=(m_{q_1},\ldots,m_{q_2})}{\underset{m'_{q_1-1} = k-1\text{ and } m'_{q_{2}+1}=l}{m_{q_1},\ldots, m_{q_2}\in [k,l]}}.
\end{array}\right. 
$$

Consider the process $\D_{n,\vec m}$ defined by
$$
\D_{n,\vec m}(s,\vec u) = \sqrt{n}\lambda_n(0,s)\lambda_n(s,1)\{C_{1:\ip{ns},\vec m}(\vec u) - C_{\ip{ns}+1:n,\vec m}(\vec u)\},\quad (s,\vec u)\in[0,1]^{d+1}.
$$

\begin{prop}\label{prop2}
Let $\vec X_1,\ldots,\vec X_n$ be $d$-dimensional random vectors with copula $C$, such that for $0\leq b_1\leq\ldots\leq b_p\leq 1$ known and $m_j=\ip{nb_j}$ for $j=1,\ldots,p$, the random vectors $\vec X_{m_{j-1}+1},\ldots, \vec X_{m_{j}}$ have m.c.d.f.s $F_{1j},\ldots,F_{dj}$. \\


Then, under Condition \ref{Cj1}, the process $\D_{n,\vec m}$ converges weakly in $\ell^{\infty}([0,1]^{d+1})$, to $\D_C$ defined in \eqref{DC}.
\end{prop}

The proof of Proposition \ref{prop2} is similar to the the proof of Proposition \ref{prop1} in which the supremas are broken on $C^2_{p+2}$ supremas.

For $s\in[0,1]$, denote $q$ the integer such that $0 =b_0\leq b_1 <\ldots<b_q< s \leq b_{q+1}<\ldots<b_p\leq b_{p+1}=1$. $\D_{n,\vec m}$ rewrites for $(s,\vec u)\in[0,1]^{d+1}$ as 
$$
\D_{n,\vec m}(s,\vec u) = \lambda_n(s,1)\{\sum_{j=1}^{q}\C_{n,m}(b_{j-1},b_j,\vec u) + \C_{n,m}(b_q,s,\vec u)\}- \lambda_n(0,s)\{\C_{n,m}(s,b_{q+1},\vec u) + \sum_{j=q+1}^{p}\C_{n,m}(b_{j},b_{j+1},\vec u)\}.
$$
For $\beta=1,\ldots,B$, consider the re-sampling versions 
$$
\check{\D}^{(\beta)}_{n,\vec m}(s,\vec u) = \lambda_n(s,1)\{\sum_{j=1}^{q}\check{\C}^{(\beta)}_{n,m}(b_{j-1},b_j,\vec u) + \check{\C}^{(\beta)}_{n,m}(b_q,s,\vec u)\}- \lambda_n(0,s)\{\check{\C}^{(\beta)}_{n,m}(s,b_{q+1},\vec u) + \sum_{j=q+1}^{p}\check{\C}^{(\beta)}_{n,m}(b_{j},b_{j+1},\vec u)\}.
$$
\begin{prop}
Under the same condition as Proposition 3.1, the conclusion of proposition 2.2 holds with $\D_{n,\vec m}$ instead of $\D_{n,m}$ and for $\beta=1,\ldots,B$, $\check{\D}^{(\beta)}_{n,\vec m}$ instead of $\check{\D}^{(\beta)}_{n,m}$.
\end{prop}

This Proposition can be proved in the same way as the prove of Proposition 2.2.

\section{Monte Carlo simulations}

In all the simulations, $d$-dimensional observations were considered,  with either a Clayton (Cl) copula or a Gumbel--Hougaard (GH) copula defined for $\vec u= (u_1,\ldots,u_d)\in[0,1]^d$ by
\begin{align*}
C_{\theta}^{Cl}(\vec u) &= \max\left(\sum_{j=1}^du_j - d+1,0\right)^{-1/\theta},\quad \theta\geq 1,\\
C_{\theta}^{GH}(\vec u) &= \exp\left(-\left[\sum_{j=1}^d\left\{-\log(u_j)\right\}^\theta\right]^{1/\theta}\right),\quad \theta>0.
\end{align*}

In an equivalent way, the Kendall's tau of bivariate margins were specified instead of the parameter $\theta$ of the copula. The Monte Carlo experiments were generated using the  \citet[R statistical system][]{Rsystem} and the \textsl{copula} package of \cite{copula} to sample the Clayton and Gumbel--Hougaard Copulas. The reader may request the corresponding routine by contacting the author. \\

In a first situation corresponding to the simulations appearing in Table~\ref{tableH0}, independent samples of sizes $n=\{50, 100, 200\}$, and dimensions $d=\{2,3\}$  are considered, where the first $\bstar=\ip{nb}$ observations, $b=\{0.1,0.25,0.5\}$,  have for marginal distributions, normal distributions $N(2,1)$ and for copula a Clayton copula or a Gumbel--Hougaard copula. The bivariate margins have a Kendall's tau of $\tau=\{0.25,0.5,0.75\}$. The last $n-\ip{nb}$ observations have for marginal distributions, normal distributions $N(0,1)$ and for copula a Clayton copula or Gumbel--Hougaard copula and the bivariate margins have a Kendall's tau of $\tau$. Typically, the samples were generated  under $H^m_0 = H_{0,c}\cap H_{1,m}$ where $H_{0,c}$ is defined in \eqref{H0c} and $H_{1,m}$ in \eqref{H1mb}. The percentage of rejection of $H^m_0$, based on $S_{n,\bstar}$ defined in~\eqref{Snb} with a level $\alpha=5\%$ were studied.\\

\begin{table}[h!]
\centering
\caption{Percentage of rejection of $H^m_0$ based on the $S_{n,\bstar}$ statistic, computed from 1000 random samples of size $n = \{50,100,200\}$ generated under $H^m_0$, where $C$ is either a $d$-dimensional Clayton (Cl) or a $d$-dimensional Gumbel-Hougaard (GH) copula with Kendall's tau of $\tau$. The $\ip{nb}$ first observations have for marginal distributions, normal distributions $N(2,1)$ and the $n-\ip{nb}$ last observations have for marginal distributions, normal distributions $N(0,1)$. The test is based on $S_n^{\bstar}$ and replications are computed using independent multipliers } 
\label{tableH0}
\begin{tabular}{rrr|rrr|rrr}
  \hline
  \multicolumn{3}{c}{} & \multicolumn{3}{c}{Cl} & \multicolumn{3}{c}{GH} \\  $n$ & $d$ & $\tau$ & $b=0.1$ & $b=0.25$ & $b=0.5$ & $b=0.1$ & $b=0.25$ & $b=0.5$  \\ \hline
50 & 2 & 0.25 & 11.3 & 9.5 & 7.4 & 7.9 & 6.9 & 7.0 \\ 
   &  & 0.50 & 16.3 & 10.4 & 4.8 & 10.4 & 7.5 & 4.0 \\ 
   &  & 0.75 & 36.1 & 11.3 & 5.4 & 18.2 & 5.8 & 1.7 \\ 
   & 3 & 0.25 & 9.3 & 9.3 & 8.1 & 5.3 & 6.3 & 5.1 \\ 
   &  & 0.50 & 10.8 & 9.3 & 7.0 & 3.9 & 3.0 & 2.5 \\ 
   &  & 0.75 & 3.4 & 2.1 & 1.7 & 1.0 & 1.3 & 0.6 \\ 
	\hline
	\hline
  100 & 2 & 0.25 & 8.2 & 7.6 & 6.6 & 4.3 & 2.9 & 3.8 \\ 
   &  & 0.50 & 8.0 & 7.0 & 4.9 & 6.6 & 5.0 & 4.4 \\ 
   &  & 0.75 & 11.5 & 6.0 & 2.1 & 5.8 & 4.1 & 1.5 \\ 
   & 3 & 0.25 & 5.7 & 6.0 & 6.3 & 5.4 & 5.5 & 5.1 \\ 
   &  & 0.50 & 8.1 & 8.5 & 7.9 & 3.6 & 3.2 & 3.4 \\ 
   &  & 0.75 & 3.8 & 2.7 & 1.1 & 1.0 & 1.2 & 0.3 \\ 
	\hline
	\hline
  200 & 2 & 0.25 & 4.7 & 5.2 & 4.5 & 4.9 & 5.1 & 5.2 \\ 
   &  & 0.50 & 4.3 & 5.5 & 4.2 & 5.3 & 5.1 & 4.4 \\ 
   &  & 0.75 & 5.1 & 4.0 & 1.8 & 4.5 & 3.1 & 1.1 \\ 
   & 3 & 0.25 & 5.2 & 5.6 & 4.8 & 4.0 & 4.1 & 4.2 \\ 
   &  & 0.50 & 7.0 & 6.7 & 5.5 & 3.7 & 3.7 & 3.1 \\ 
   &  & 0.75 & 3.6 & 3.2 & 2.8 & 2.2 & 2.0 & 1.1 \\ 
   \hline
\end{tabular}
\end{table}

These percentages of rejection are appreciably closed around $\alpha=5\%$, except for the cases with $b=0.1$ and $n=50,100$ where estimations are calculated on 5 and 10 observations. This may explain the too high percentage of rejection of $H^m_0$. \\

In Table~\ref{tableH1Cl} and~\ref{tableH1GH}, independent samples of sizes $n=\{50,100,200\}$, and dimensions $d=\{2,3\}$ are considered where the first $\bstar=\ip{nb}$ observations, $b=\{0.1,0.25,0.5\}$ have for marginal distributions, normal distributions $N(2,1)$ and for copula a Clayton copula or Gumbel--Hougaard copula. The bivariate margins have a Kendall's tau of $\tau=0.2$. The last $n-\ip{nb}$ observations have for marginal distributions, normal distributions $N(0,1)$ and for copula Clayton copula or Gumbel--Hougaard copula where the bivariate margins have a Kendall's tau of $\tau=\{0.4,0.6\}$. Typically, the samples were generated under alternative hypotheses $H^{\bstar}_A=H_{1,c}\cap H_{1,m}$ where $H_{1,m}$ is defined in~(\ref{H1mb}) and with
\begin{multline}\label{H1c}
H_{1,c}:\quad \text{There exist }k\in\{1,\ldots,n-1\} \text{ and two copulas } C_1 \text{ and }C_2,\\ \text{ such that }\vec X_1,\ldots,\vec X_k \text{ have copula }C_1 \text{ and }\vec X_{k+1},\ldots,\vec X_n \text{ have copula }C_2.
\end{multline}

\setlength{\footskip}{120pt} 
\begin{table}
\centering
\caption{Percentage of rejection of $H^m_0$  computed from 1000 random samples of size $n = \{50,100,200\}$ generated under $H_A = H_{1,c}\cap H_{1,m}$ defined in \eqref{H1c} and \eqref{H1mb}, where the first $\ip{nt}$ observations, $t\in\{0.1,0.25,05\}$ have for copula a $d$-dimensional Clayton (Cl) copula with Kendall's tau of $0.2$ and the last $n-\ip{nt}$ have for copula a $d$-dimensional Clayton copula (Cl) with Kendall's tau of $\tau$. The first $\bstar=\ip{nb}$ observations, $b\in\{0.1,0.25,0.5\}$ have for marginal distributions, normal distributions $N(0,1)$ and the $n-\bstar$ last observations have for marginal distributions, normal distributions $N(2,1)$. Two different tests (based on $S_n$ and $S_n^{\bstar}$) are compared and replications are computed using independent multipliers  } 
\label{tableH1Cl}
\begin{tabular}{rrrr||rrrrrr}
  \hline
   \multicolumn{4}{c}{} & \multicolumn{1}{c}{$S_{n,\bstar}$} & \multicolumn{1}{c}{$S_n$}& \multicolumn{1}{c}{$S_{n,\bstar}$} & \multicolumn{1}{c}{$S_n$}& \multicolumn{1}{c}{$S_{n,\bstar}$} & \multicolumn{1}{c}{$S_n$} \\ $n$ & $d$ & $t$ & $\tau$ & $b=0.1$ & $b=0.1$ & $b=0.25$ & $b=0.25$ & $b=0.5$ & $b=0.5$  \\ \hline
50 & 2 & 0.10 & 0.4 & 15.9 & 10.5 & 10.5 & 15.1 & 7.7 & 6.7 \\ 
   &  &  & 0.6 & 39.9 & 28.5 & 23.0 & 45.0 & 10.7 & 8.4 \\ 
   &  & 0.25 & 0.4 & 21.6 & 19.3 & 18.0 & 21.1 & 11.6 & 9.6 \\ 
   &  &  & 0.6 & 56.6 & 48.5 & 48.5 & 60.0 & 31.8 & 19.9 \\ 
   &  & 0.50 & 0.4 & 23.8 & 23.8 & 22.0 & 25.9 & 17.0 & 15.3 \\ 
   &  &  & 0.6 & 59.0 & 56.3 & 55.7 & 63.7 & 46.3 & 43.2 \\ 
   & 3 & 0.10 & 0.4 & 12.8 & 12.6 & 11.3 & 15.9 & 8.5 & 6.0 \\ 
   &  &  & 0.6 & 26.0 & 25.4 & 20.7 & 40.9 & 13.6 & 8.7 \\ 
   &  & 0.25 & 0.4 & 20.0 & 21.0 & 21.7 & 25.0 & 14.8 & 10.5 \\ 
   &  &  & 0.6 & 57.8 & 62.6 & 56.2 & 67.6 & 47.2 & 30.3 \\ 
   &  & 0.50 & 0.4 & 28.7 & 31.2 & 27.5 & 35.6 & 24.9 & 24.6 \\ 
   &  &  & 0.6 & 75.9 & 79.4 & 76.9 & 82.3 & 73.5 & 72.3 \\
	\hline
	\hline
  100 & 2 & 0.10 & 0.4 & 10.5 & 17.5 & 9.5 & 27.4 & 6.4 & 8.1 \\ 
   &  &  & 0.6 & 35.4 & 47.9 & 28.4 & 72.8 & 19.5 & 25.4 \\ 
   &  & 0.25 & 0.4 & 23.8 & 33.1 & 24.0 & 40.6 & 17.3 & 13.6 \\ 
   &  &  & 0.6 & 78.0 & 84.2 & 74.7 & 90.1 & 62.7 & 47.8 \\ 
   &  & 0.50 & 0.4 & 32.4 & 38.3 & 32.4 & 45.5 & 29.4 & 30.0 \\ 
   &  &  & 0.6 & 84.0 & 89.7 & 85.1 & 95.2 & 81.9 & 85.0 \\ 
   & 3 & 0.10 & 0.4 & 11.8 & 16.3 & 11.5 & 22.9 & 9.7 & 7.9 \\ 
   &  &  & 0.6 & 36.6 & 49.4 & 33.1 & 69.0 & 25.7 & 21.8 \\ 
   &  & 0.25 & 0.4 & 30.8 & 40.0 & 31.1 & 45.1 & 26.9 & 19.7 \\ 
   &  &  & 0.6 & 86.7 & 91.9 & 87.0 & 94.4 & 82.4 & 63.6 \\ 
   &  & 0.50 & 0.4 & 43.1 & 53.0 & 42.4 & 59.1 & 41.9 & 41.8 \\ 
   &  &  & 0.6 & 95.8 & 97.0 & 95.4 & 98.5 & 95.6 & 96.1 \\ 
	\hline
	\hline
  200 & 2 & 0.10 & 0.4 & 11.2 & 30.2 & 11.7 & 45.8 & 10.2 & 16.6\\
   &  &  & 0.6 & 49.0 & 81.4 & 42.5 & 96.1 & 37.7 & 70.1 \\ 
   &  & 0.25 & 0.4 & 35.6 & 56.3 & 36.2 & 68.7 & 31.1 & 31.7 \\ 
   &  &  & 0.6 & 93.4 & 98.0 & 93.5 & 99.7 & 91.1 & 89.1 \\ 
   &  & 0.50 & 0.4 & 50.3 & 65.6 & 49.8 & 76.0 & 49.3 & 58.4 \\ 
   &  &  & 0.6 & 99.1 & 99.9 & 98.9 & 100.0 & 99.0 & 99.5 \\ 
   & 3 & 0.10 & 0.4 & 12.7 & 29.2 & 12.6 & 44.1 & 11.2 & 14.3 \\ 
   &  &  & 0.6 & 61.5 & 81.6 & 59.7 & 94.4 & 57.9 & 59.8 \\ 
   &  & 0.25 & 0.4 & 46.2 & 66.5 & 48.8 & 72.7 & 45.0 & 36.7 \\ 
   &  &  & 0.6 & 98.7 & 99.9 & 99.2 & 99.9 & 99.0 & 94.6 \\ 
   &  & 0.50 & 0.4 & 69.0 & 82.3 & 69.3 & 88.8 & 70.2 & 75.2 \\ 
   &  &  & 0.6 & 100.0 & 100.0 & 100.0 & 100.0 & 100.0 & 100.0 \\ 
   \hline
\end{tabular}
\end{table}
\begin{table}
\centering
\caption{Percentage of rejection of $H^m_0$  computed from 1000 random samples of size $n = \{50,100,200\}$ generated under $H_A = H_{1,c}\cap H_{1,m}$ defined in \eqref{H1c} and \eqref{H1mb}, where the first $\ip{nt}$ observations, $t\in\{0.1,0.25,05\}$ have for copula a $d$-dimensional Gumbel--Hougaard (GH) copula with Kendall's tau of $0.2$ and the last $n-\ip{nt}$ have for copula a $d$-dimensional Gumbel--Hougaard (GH) with Kendall's tau of $\tau$. The  first $\bstar = \ip{nb}$ observations, $b\in\{0.1,0.25,0.5\}$ have for marginal distributions, normal distributions $N(0,1)$ and the $n-\bstar$ last observations have for marginal distributions, normal distributions $N(2,1)$. Two different tests (based on $S_n$ and $S_n^{\bstar}$) are compared and replications are computed using independent multipliers  } 
\label{tableH1GH}
\begin{tabular}{rrrr||rrrrrr}
  \hline
   \multicolumn{4}{c}{} & \multicolumn{1}{c}{$S_{n,\bstar}$} & \multicolumn{1}{c}{$S_n$}& \multicolumn{1}{c}{$S_{n,\bstar}$} & \multicolumn{1}{c}{$S_n$}& \multicolumn{1}{c}{$S_{n,\bstar}$} & \multicolumn{1}{c}{$S_n$} \\ $n$ & $d$ & $t$ & $\tau$ & $b=0.1$ & $b=0.1$ & $b=0.25$ & $b=0.25$ & $b=0.5$ & $b=0.5$  \\ \hline
50 & 2 & 0.10 & 0.4 & 11.2 & 9.0 & 8.5 & 13.7 & 5.6 & 4.8 \\ 
   &  &  & 0.6 & 23.5 & 22.6 & 15.1 & 31.4 & 8.9 & 6.5 \\ 
   &  & 0.25 & 0.4 & 18.2 & 14.9 & 15.2 & 19.5 & 11.1 & 7.3 \\ 
   &  &  & 0.6 & 45.9 & 41.8 & 39.4 & 49.1 & 28.8 & 15.8 \\ 
   &  & 0.50 & 0.4 & 20.7 & 17.9 & 19.9 & 22.0 & 17.1 & 15.7 \\ 
   &  &  & 0.6 & 52.9 & 55.0 & 49.5 & 60.1 & 46.5 & 43.3 \\ 
   & 3 & 0.10 & 0.4 & 7.1 & 8.2 & 6.8 & 12.2 & 6.1 & 5.8 \\ 
   &  &  & 0.6 & 13.4 & 16.3 & 10.3 & 26.1 & 8.2 & 4.7 \\ 
   &  & 0.25 & 0.4 & 12.2 & 15.8 & 14.0 & 19.0 & 11.5 & 8.5 \\ 
   &  &  & 0.6 & 41.3 & 49.6 & 41.6 & 54.5 & 36.1 & 22.9 \\ 
   &  & 0.50 & 0.4 & 20.9 & 25.9 & 21.5 & 27.9 & 21.8 & 18.9 \\ 
   &  &  & 0.6 & 66.1 & 72.1 & 68.0 & 73.8 & 63.9 & 62.4 \\ 
	\hline
	\hline
  100 & 2 & 0.10 & 0.4 & 9.1 & 11.5 & 7.9 & 21.6 & 6.6 & 6.6 \\ 
   &  &  & 0.6 & 26.5 & 44.0 & 18.3 & 59.5 & 14.1 & 17.2 \\ 
   &  & 0.25 & 0.4 & 19.4 & 28.3 & 20.5 & 33.1 & 17.1 & 11.6 \\ 
   &  &  & 0.6 & 66.7 & 76.0 & 62.4 & 84.2 & 55.9 & 39.3 \\ 
   &  & 0.50 & 0.4 & 26.7 & 34.6 & 27.2 & 41.5 & 24.8 & 25.8 \\ 
   &  &  & 0.6 & 81.6 & 87.5 & 80.2 & 91.0 & 79.3 & 81.2 \\ 
   & 3 & 0.10 & 0.4 & 6.9 & 11.9 & 6.6 & 18.0 & 6.2 & 5.3 \\ 
   &  &  & 0.6 & 25.1 & 38.8 & 21.8 & 55.1 & 19.5 & 12.6 \\ 
   &  & 0.25 & 0.4 & 21.5 & 29.4 & 22.8 & 34.0 & 21.2 & 14.5 \\ 
   &  &  & 0.6 & 79.1 & 87.6 & 79.6 & 90.0 & 76.8 & 55.1 \\ 
   &  & 0.50 & 0.4 & 39.3 & 50.8 & 38.5 & 55.1 & 39.3 & 39.7 \\ 
   &  &  & 0.6 & 93.2 & 95.9 & 93.4 & 96.9 & 92.8 & 91.9 \\ 
	\hline
	\hline
  200 & 2 & 0.10 & 0.4 & 8.7 & 26.4 & 8.2 & 42.4 & 7.8 & 15.7 \\ 
   &  &  & 0.6 & 40.5 & 79.3 & 35.2 & 92.7 & 33.4 & 56.2 \\ 
   &  & 0.25 & 0.4 & 31.3 & 52.1 & 29.2 & 62.6 & 28.3 & 26.2 \\ 
   &  &  & 0.6 & 91.1 & 98.3 & 91.5 & 99.8 & 90.3 & 80.4 \\ 
   &  & 0.50 & 0.4 & 46.3 & 63.0 & 45.2 & 72.6 & 44.2 & 51.7 \\ 
   &  &  & 0.6 & 99.1 & 99.8 & 99.0 & 100.0 & 99.4 & 99.6 \\ 
   & 3 & 0.10 & 0.4 & 8.2 & 25.3 & 9.9 & 34.9 & 9.6 & 9.5 \\ 
   &  &  & 0.6 & 52.5 & 80.4 & 48.6 & 89.4 & 48.6 & 40.0 \\ 
   &  & 0.25 & 0.4 & 40.2 & 62.2 & 42.0 & 64.0 & 42.0 & 31.0 \\ 
   &  &  & 0.6 & 98.0 & 99.4 & 98.1 & 99.6 & 97.9 & 91.1 \\ 
   &  & 0.50 & 0.4 & 65.9 & 79.5 & 65.5 & 85.5 & 66.1 & 69.7 \\ 
   &  &  & 0.6 & 99.9 & 100.0 & 100.0 & 100.0 & 100.0 & 100.0 \\ 
   \hline
\end{tabular}
\end{table}

The break times $k=\ip{nt}$, $t\in\{0.1,0.25,0.5\}$ are considered. The percentages of rejection of the hypothesis $H^m_0$ with a level $\alpha=5\%$ are studied. In the same way, the test for $H_0 = H_{0,c}\cap H_{0,m}$ where $H_{0,m}$ is defined in~(\ref{H0m}) is considered, based on $S_n$ described in equation (\ref{Sn}).
The percentages of rejection of $H^m_0$ based on $S_{n,\bstar}$ are closed to the percentages of rejection of $H_0$ based on $S_n$. More exactly the percentages of rejection of $H^m_0$ based on $S_{n,\bstar}$ are generally smaller than percentage of rejection of $H_0$ based $S_n$ for $b=\{0.1,0.25\}$ and larger for $b=0.5$ and $t\in\{0.25,0.5\}$. \\

Recall that with the hypothesis of a break time $\bstar$ known in the m.c.d.f.s, the rejection of $H_0$ using $S_n$ does not allow for a conclusion of a break in the copula of observations contrary to the rejection of $H^m_0$ using $S_{n,\bstar}$.\\

\section{Discussions and specific situation}
\subsection{A Strong mixing condition}
Suppose that the random vectors $\vec X_1,\ldots,\vec X_n$ are drawn from sequences of weakly dependent vectors, in the sense of $\alpha$-mixing dependence (strong mixing dependence) introduced in \cite{Ros56}:
\begin{definition}
Let $(\vec X_i)_{i\in\Z}$ a sequence of random vectors, and for $a,b\in\bar \Z = \Z\cup\{\pm \infty\}$, denote by $\mathcal F_a^b$ the $\sigma$-field generated by $(\vec X_i)_{a\leq i\leq b}$. The sequence of $\alpha$-mixing coefficients $(\alpha_r)_{r\in\mathbb N}$ is defined by
$$
\alpha_r = \sup_{t\in\Z}\sup_{A\in\mathcal F_{-\infty}^t, B\in\mathcal F_{t+r}^{+\infty}}\left|P(A\cap B) - P(A)P(B)\right|, \quad r\in\mathbb N.
$$
The sequence $(\vec X_i)_{i\in\Z}$ will be said to be strongly mixing as soon as $\alpha_r\underset{r\to +\infty}{\longrightarrow} 0$.
\end{definition}

\begin{enumerate}
\item The propositions~\ref{prop1} and ~\ref{prop2} remain true if we suppose that the marginal probability integral transforms $\vec U_{1,m},\ldots,  \vec U_{n,m}$ defined in \eqref{Uim} are drawn from a strictly stationary sequence $(\vec U_{i})_{i\in\Z}$ whose strong mixing coefficients satisfy $\alpha_r=O(n^{-a})$, $a>1$. In this case the covariance structure of $\Z_C$ is given by $\cov(\Z_C(s,\vec u),\Z_C(t,\vec v)) = \min(s,t)\sum_{k\in\Z}\cov(\1(\vec U_{0}\leq \vec u),\1(\vec U_{k}\leq \vec v))$, $(s,\vec u)$, $(t,\vec v)\in[0,1]^{d+1}$.
%
\item  The propositions~\ref{prop:Cnmboot} and 3.2  remains true if we suppose that the marginal probability integral transforms $\vec U_{1,m},\ldots,  \vec U_{n,m}$ defined in \eqref{Uim} are drawn from a strictly stationary sequence $(\vec U_{i})_{i\in\Z}$ whose strong mixing coefficients satisfy $\alpha_r=O(n^{-a})$,  $a>3+3d/2$ and we consider dependent multipliers satisfy (M1)--(M3) appearing in \citet[section 2]{BucKoj14} with $\ell_n = O(n^{1/2-\gamma})$ for some $0<\gamma<1/2$ instead of independent multipliers.

\end{enumerate}

These situations have been studied in Tables~\ref{H0dep} and \ref{HAdep}; sequences of  multipliers were simulated  using the procedure of \citep[The moving average approach, Section 6.1]{BucKoj14}. A standard normal sequence of i.i.d. random variables was used in the construction of multipliers.  The value of the bandwidth appearing in the condition (M2) was automatically selected by the procedure described in \citep[Section 5]{BucKoj14} by using the R function \textsl{bOptEmpProc} of \textsl{npcp} package (\cite{npcp}). The "combining" function $\psi$ appearing in this same procedure was arbitrarily chosen as $\psi =maximum$ \citep[see][ Section 4]{PolWhi04}. 
Finally the function $\varphi$ appearing in the condition (M3) was the convolution product $\varphi(x) = \kappa_P \star \kappa_P(2x)/\kappa_P \star \kappa_P(0)$, where $\kappa_P = (1 - 6x^2+6|x|^3)\1(|x|\leq 1/2) + 2(1-|x|)^3\1(1/2<|x|\leq 1)$, $x\in\R$.\\

In Table~\ref{H0dep} and Table \ref{HAdep}, dependent samples of sizes $n=\{100,200\}$ and $d=\{2,3\}$ are considered with a break in the variance of marginal distributions at time $m=\ip{nb}$, $b=\{0.1,0.25,0.5\}$. The samples were generated under $H^m_0 = H_{0,c}\cap H_{1,m}$ in Table~\ref{H0dep} and under  $(\neg H_{0,c})\cap H_{1,m}$ in Table \ref{HAdep} with a break in the copula at time $k=\ip{nt}$, $t=\{0.1,0.25,0.5\}$. The data are generated from two autoregressive models (AR1) defined by:
\begin{equation}\label{AR1}\tag{AR1}
X_{i+1,j} = 0.5 X_{i,j} + \varepsilon_{i+1,j},\quad j=1,\ldots,d \text{ and }i\in \Z\quad \text{with }\varepsilon_{ij}\text{ white noise.}
\end{equation}
For $i=1,\ldots m$ and $j=1,\ldots,d$ the chosen white noises are $N(0,1)$, and for $i=m+1,\ldots n$ the chosen white noises are $N(0,16)$. 
For $b\leq t$, the sample is obtained in the following way:
let $\vec U_{1},\ldots,\vec U_{k+200}$ be a d-variate i.i.d. sample from the copula $C_1$ and $\vec U_{k+201},\ldots,\vec U_{n+200}$ be a d-variate i.i.d. sample from the copula $C_2$ (under $H^m_0$, $C_1=C_2$). For $i=1,\ldots,m+100$, let $\vec \varepsilon_i = (\Phi^{-1}(U_{i1}),\ldots,\Phi^{-1}(U_{id}))$ where $\Phi$ is the c.d.f. of the standard normal distribution, and for $i=m+101,\ldots,n+200$, $\vec \varepsilon_i = (4\times\Phi^{-1}(U_{i1}),\ldots,4\times\Phi^{-1}(U_{id}))$.
Then, $\vec X_1 = \vec \varepsilon_1$, $\vec X_{m+101} = \vec \varepsilon_{m+101}$  and for $i=1,\ldots,m+99$, and $i=m+101,\ldots,n+199$ compute recursively
$$
\vec X_{i+1} = 0.5\vec X_{i} + \varepsilon_{i+1}.
$$
Finally, we remove the observations $\vec X_1$ to $\vec X_{100}$ and $\vec X_{m+101}$ to $\vec X_{m+200}$. For $b>t$, the sample is obtained in the similar way.

\begin{table}[t!]
\centering
\caption{Percentage of rejection of $H^m_0$  computed from 1000 samples of size $n = \{100,200\}$ generated under $H^m_0$ and from two (AR1) models, with $d$-dimensional Clayton (Cl) or Gumbel--Hougaard (GH) stationary copula with Kendall's tau of $\tau$. The $\bstar=\ip{nb}$ first observations, $b\in\{0.1,0.25,0.5\}$  have stationary margins $N(0,1)$ and the $n-\ip{nb}$ last observations have stationary margins $N(0,16)$. The test is based on $S_n^{\bstar}$ and replications are computed using dependent multipliers } 
\label{H0dep}
\begin{tabular}{rrr||rrr||rrr}
  \hline
  \multicolumn{3}{c}{} & \multicolumn{3}{c}{Cl} & \multicolumn{3}{c}{GH} \\  $n$ & $d$ & $\tau$ & $b=0.1$ & $b=0.25$ & $b=0.5$ & $b=0.1$ & $b=0.25$ & $b=0.5$  \\ \hline
100 & 2 & 0.25 & 9.0 & 10.6 & 9.6 & 8.3 & 10.3 & 8.3 \\ 
   &  & 0.50 & 11.2 & 8.2 & 8.5 & 11.5 & 8.1 & 6.2 \\ 
   &  & 0.75 & 12.2 & 11.2 & 5.1 & 12.9 & 8.1 & 5.0 \\ 
   & 3 & 0.25 & 9.7 & 9.9 & 10.6 & 7.0 & 10.3 & 8.8 \\ 
   &  & 0.50 & 9.9 & 10.7 & 8.1 & 5.0 & 6.8 & 7.2 \\ 
   &  & 0.75 & 4.5 & 3.4 & 3.3 & 2.5 & 2.8 & 1.4 \\ 
	\hline
	\hline
  200 & 2 & 0.25 & 5.9 & 8.0 & 6.5 & 7.6 & 6.5 & 5.6 \\ 
   &  & 0.50 & 6.1 & 7.1 & 5.1 & 3.6 & 5.4 & 4.8 \\ 
   &  & 0.75 & 4.3 & 4.6 & 3.4 & 3.4 & 2.9 & 0.7 \\ 
   & 3 & 0.25 & 6.2 & 4.8 & 6.4 & 5.0 & 6.3 & 8.2 \\ 
   &  & 0.50 & 5.2 & 5.7 & 4.5 & 4.9 & 4.3 & 4.6 \\ 
   &  & 0.75 & 2.1 & 3.6 & 1.9 & 1.0 & 1.1 & 0.6 \\ 
   \hline
\end{tabular}
\end{table}
\begin{table}[t!]
\centering
\caption{Percentage of rejection of $H^m_0$  computed from 1000 samples of size $n = \{100,200\}$  generated under $\neg H_{0,c}\cap H_{1,m}$ and from two (\ref{AR1}) models, where the first $\ip{nt}$ observations, $t\in\{0.1,0.25,05\}$ have for stationary copula a $d$-dimensional Clayton (Cl) copula (resp. Gumbel--Hougaard Copula) with Kendall's tau of $0.2$ and the last $n-\ip{nt}$ have for stationary copula a bidimensional Clayton copula (Cl) (resp. Gumbel--Hougaard Copula) with Kendall's tau of $\tau$. The $\ip{nb}$ first observations have m.c.d.f. $N(0,1)$ and the $n-\ip{nb}$ last observations have m.c.d.f. $N(0,16)$.  The test is based on $S_n^{\bstar}$ and replications are computed using dependent multipliers  } 
\label{HAdep}
\begin{tabular}{rrrr||rrr||rrr}
  \hline
  \multicolumn{4}{c}{} & \multicolumn{3}{c}{Cl} & \multicolumn{3}{c}{GH} \\  $n$ & $d$ & $t$  & $\tau$ & $b=0.1$ & $b=0.25$ & $b=0.5$ & $b=0.1$ & $b=0.25$ & $b=0.5$  \\ \hline
100 & 2 & 0.10 & 0.4 & 14.8 & 13.4 & 10.5 & 12.7 & 12.3 & 9.1 \\ 
   &  &  & 0.6 & 36.5 & 23.9 & 10.7 & 33.1 & 19.6 & 9.4 \\ 
   &  & 0.25 & 0.4 & 24.6 & 25.4 & 17.5 & 22.8 & 20.4 & 15.4 \\ 
   &  &  & 0.6 & 62.5 & 62.8 & 45.9 & 62.0 & 55.7 & 45.1 \\ 
   &  & 0.50 & 0.4 & 31.9 & 32.2 & 25.2 & 27.0 & 24.8 & 25.8 \\ 
   &  &  & 0.6 & 74.7 & 71.0 & 69.1 & 73.8 & 72.1 & 66.4 \\ 
   & 3 & 0.10 & 0.4 & 12.2 & 13.2 & 11.8 & 11.3 & 10.6 & 9.0 \\ 
   &  &  & 0.6 & 34.9 & 23.1 & 13.8 & 26.0 & 18.3 & 11.1 \\ 
   &  & 0.25 & 0.4 & 23.7 & 31.5 & 20.9 & 23.7 & 25.3 & 18.7 \\ 
   &  &  & 0.6 & 74.0 & 77.0 & 64.9 & 67.3 & 72.6 & 59.0 \\ 
   &  & 0.50 & 0.4 & 38.8 & 40.8 & 36.1 & 32.3 & 35.0 & 35.5 \\ 
   &  &  & 0.6 & 85.1 & 87.6 & 86.8 & 85.6 & 86.6 & 84.8 \\ 
	\hline
	\hline
  200 & 2 & 0.10 & 0.4 & 10.8 & 10.4 & 9.6 & 10.5 & 8.6 & 7.9 \\ 
   &  &  & 0.6 & 39.9 & 21.9 & 11.3 & 36.3 & 20.1 & 12.2 \\ 
   &  & 0.25 & 0.4 & 27.9 & 29.5 & 22.1 & 24.2 & 26.6 & 20.6 \\ 
   &  &  & 0.6 & 80.9 & 81.5 & 72.1 & 79.8 & 82.0 & 71.9 \\ 
   &  & 0.50 & 0.4 & 39.5 & 37.5 & 35.7 & 35.7 & 38.4 & 32.6 \\ 
   &  &  & 0.6 & 91.9 & 92.5 & 90.3 & 92.5 & 90.9 & 92.0 \\ 
   & 3 & 0.10 & 0.4 & 14.9 & 11.4 & 8.9 & 10.4 & 9.6 & 8.4 \\ 
   &  &  & 0.6 & 45.8 & 26.6 & 17.6 & 45.2 & 21.4 & 15.6 \\ 
   &  & 0.25 & 0.4 & 37.3 & 38.9 & 31.9 & 30.6 & 35.0 & 29.1 \\ 
   &  &  & 0.6 & 93.5 & 94.1 & 86.8 & 90.7 & 93.6 & 86.9 \\ 
   &  & 0.50 & 0.4 & 54.6 & 53.0 & 54.0 & 51.5 & 51.8 & 51.9 \\ 
   &  &  & 0.6 & 99.3 & 98.0 & 98.0 & 98.5 & 98.9 & 98.4  \\ 
   \hline
\end{tabular}
\end{table}

From $n=200$, it can seen in Table~\ref{H0dep} that the percentages of rejection of $H^m_0$ are appreciably closed around $\alpha=5\%$ whereas for the whole of break scenarios in copula (Table~\ref{HAdep}), the percentages of rejection of $H^m_0$ are relatively high. 

\subsection{Specific situation}

As an illustration, the bivariate log-returns computed from closing daily quotes of the Dow Jones Industrial Average and the Nasdaq Composite for the years 1987 and 1988 have been studied. This is an interesting situation because the data highlight a change in the m.c.d.f.s at time $m=202$ (1987-10-19, corresponding to the "Black Monday"). A Cram\'er--von Mises test \citep[see for example][]{HolKojQue13} can allow to confirm this change. Using the procedure described in Section~2, an approximate p-value of 0,201 was obtained and no evidence against $H_{0,c}$ is reported. 

Because a marginal gradual change or a multiple marginal change could lead to a rejection of $H^m_0$, in the case of rejection of $H^m_{0}$, the hypothesis of a unique change in marginal distribution should be confirmed.


\subsection{Case of unknown marginal break}

It seem interesting not to fix the break time $m=\ip{nb}$ and aggregating this term in the best possible way. Nevertheless,  such aggregation would be expensive in simulation time, thus the performance of the algorithm should be improved before. 

Another interesting way of a future research will be to consider an estimation of the unknown break time instead of $m$ in $S_{n,m}$ and study the associated statistic.
 
\newpage
\textbf{References}

\begin{thebibliography}{16}
\providecommand{\natexlab}[1]{#1}
\providecommand{\url}[1]{\texttt{#1}}
\expandafter\ifx\csname urlstyle\endcsname\relax
  \providecommand{\doi}[1]{doi: #1}\else
  \providecommand{\doi}{doi: \begingroup \urlstyle{rm}\Url}\fi

\bibitem[Bai(1994)]{Bai94}
J.~Bai.
\newblock Weak convergence of the sequential empirical processes of residuals
  in {ARMA} models.
\newblock \emph{The Annals of Statistics}, 22\penalty0 (4):\penalty0
  2051--2061, 1994.

\bibitem[B\"ucher and Kojadinovic(2015)]{BucKoj14}
A.~B\"ucher and I.~Kojadinovic.
\newblock A dependent multiplier bootstrap for the sequential empirical copula
  process under strong mixing.
\newblock \emph{Bernoulli}, 2015.

\bibitem[{B\"ucher} et~al.(2014){B\"ucher}, {Kojadinovic}, {Rohmer}, and
  {Segers}]{BucKojRohSeg14}
A.~{B\"ucher}, I.~{Kojadinovic}, T.~{Rohmer}, and J.~{Segers}.
\newblock Detecting changes in cross-sectional dependence in multivariate time
  series.
\newblock \emph{Journal of Multivariate Analysis}, 132\penalty0 (0):\penalty0
  111 -- 128, 2014.

\bibitem[Cs\"org\H{o} and Horv{\'a}th(1997)]{CsoHor97}
M.~Cs\"org\H{o} and L.~Horv{\'a}th.
\newblock \emph{Limit theorems in change-point analysis}.
\newblock Wiley Series in Probability and Statistics. John Wiley \& Sons,
  Chichester, UK, 1997.

\bibitem[Deheuvels(1979)]{Deh79}
P.~Deheuvels.
\newblock La fonction de d\'ependance empirique et ses propri\'et\'es: un test
  non param\'etrique d'ind\'ependance.
\newblock \emph{Acad. Roy. Belg. Bull. Cl. Sci. 5th Ser.}, 65:\penalty0
  274--292, 1979.

\bibitem[Hofert et~al.(2015)Hofert, Kojadinovic, M\"achler, and Yan]{copula}
M.~Hofert, I.~Kojadinovic, M.~M\"achler, and J.~Yan.
\newblock \emph{copula: {M}ultivariate dependence with copulas}, 2015.
\newblock URL \url{http://CRAN.R-project.org/package=copula}.
\newblock {R} package version 0.999-13.

\bibitem[Holmes et~al.(2013)Holmes, Kojadinovic, and Quessy]{HolKojQue13}
M.~Holmes, I.~Kojadinovic, and J-F. Quessy.
\newblock Nonparametric tests for change-point detection \`a la {G}ombay and
  {H}orv\'ath.
\newblock \emph{Journal of Multivariate Analysis}, 115:\penalty0 16--32, 2013.

\bibitem[Inoue(2001)]{Ino01}
A.~Inoue.
\newblock Testing for distributional change in time series.
\newblock \emph{Econometric Theory}, 17\penalty0 (1):\penalty0 156--187, 2001.

\bibitem[Kojadinovic(2014)]{npcp}
I.~Kojadinovic.
\newblock \emph{npcp: Some nonparametric tests for change-point detection in
  (multivariate) observations}, 2014.
\newblock {R} package version 0.1-1.

\bibitem[Kojadinovic et~al.(2011)Kojadinovic, Segers, and Yan]{KojSegYan11}
I.~Kojadinovic, J.~Segers, and J.~Yan.
\newblock Large-sample tests of extreme-value dependence for multivariate
  copulas.
\newblock \emph{The Canadian Journal of Statistics}, 39\penalty0 (4):\penalty0
  703--720, 2011.

\bibitem[Politis and White(2004)]{PolWhi04}
D.N. Politis and H.~White.
\newblock Automatic block-length selection for the dependent bootstrap.
\newblock \emph{Econometric Reviews}, 23\penalty0 (1):\penalty0 53--70, 2004.

\bibitem[{R Development Core Team}(2013)]{Rsystem}
{R Development Core Team}.
\newblock \emph{{R}: {A} Language and Environment for Statistical Computing}.
\newblock R Foundation for Statistical Computing, Vienna, Austria, 2013.
\newblock URL \url{http://www.R-project.org}.
\newblock {ISBN} 3-900051-07-0.

\bibitem[Rosenblatt(1956)]{Ros56}
Murray Rosenblatt.
\newblock A central limit theorem and a strong mixing condition.
\newblock \emph{Proceedings of the National Academy of Sciences of the United
  States of America}, 42\penalty0 (1):\penalty0 43, 1956.

\bibitem[Segers(2012)]{Seg12}
J.~Segers.
\newblock Asymptotics of empirical copula processes under nonrestrictive
  smoothness assumptions.
\newblock \emph{Bernoulli}, 18:\penalty0 764--782, 2012.

\bibitem[Sklar(1959)]{Skl59}
A.~Sklar.
\newblock Fonctions de r\'epartition \`a $n$ dimensions et leurs marges.
\newblock \emph{Publications de l'Institut de Statistique de l'Universit\'e de
  Paris}, 8:\penalty0 229--231, 1959.

\bibitem[{van der Vaart} and Wellner(2000)]{vanWel96}
A.W. {van der Vaart} and J.A. Wellner.
\newblock \emph{Weak convergence and empirical processes}.
\newblock Springer, New York, 2000.
\newblock Second edition.

\end{thebibliography}

\newpage
\appendix

\begin{appendix}
    \setcounter{lemma}{0}
    \renewcommand{\thelemma}{\Alph{section}\arabic{lemma}}
\section{Proof of Proposition \ref{prop1} and Proposition \ref{prop:Cnmboot}}
\newproof{pop}{Proof of Proposition \ref{prop1}}
\begin{pop}
In the sequel, for $k>l$ the empirical copulas $C_{k:l}$ and $C_{k:l,m}$ are considered as null by convention. Let $b\in(0,1)$ such that $\bstar=\ip{nb}$.

Let us consider the two-sided sequential empirical copula process defined in \citet{BucKoj14,BucKojRohSeg14} by:
\begin{equation}\label{Cn}
\C_n(s,t,\vec u) = \sqrt{n}\lambda_n(s,t)\{C_{\ip{ns}+1:\ip{nt}}(\vec u) - C(\vec u)\},\quad (s,t,\vec u)\in\Delta\times[0,1]^d,
\end{equation}
where $\Delta = \{(s,t)\in[0,1]^2|s\leq t\}$. If $\vec X_1,\ldots,\vec X_n$ are drawn from a i.i.d. sequence $(\vec X_i)_{i\in\Z}$ with continuous margins $F_1,\ldots,F_d$ then (Proposition~3.3 of \cite{BucKojRohSeg14}), 
\begin{equation}\label{buckojseg}
\sup_{(s,t,\vec u)\in\Delta\times[0,1]^d}|\C_{n}(s,t,\vec u) - \tilde\C_{n}(s,t,\vec u)|\p 0,
\end{equation}

with 
\begin{equation}\label{tildeCn}
\tilde\C_{n}(s,t,\vec u) = \B_n(s,t,\vec u)-\sum_{j=1}^d\dot C_j(\vec u)\B_n(s,t,\vec u^{\{j\}}),\qquad (s,t,\vec u)\in\Delta\times[0,1]^d,
\end{equation}
and
\begin{equation}\label{Bn}
\B_n(s,t,\vec u) = \frac{1}{\sqrt{n}}\sum_{i=\ip{ns}+1}^{\ip{nt}}\{\1(\vec U_i\leq \vec u) - C( \vec u)\}, \qquad (s,t,\vec u)\in\Delta\times[0,1]^d,
\end{equation}
where the vector $\vec U_i$, $i=1,\ldots,n$, are the vectors of the probability integral tranforms $\vec U_i = (F_1(X_{i1}),\ldots,F_d(X_{id}))$.\\

Here, we only suppose that $\vec X_1,\ldots,\vec X_n$ have same copula $C$.

%
\noindent Define the process $(s,t,\vec u)\mapsto\tilde\C_{n,m}(s,t,\vec u)$ similarly to $(s,t,\vec u)\mapsto\tilde\C_{n}(s,t,\vec u)$ in \eqref{tildeCn} and \eqref{Bn} using the vectors $\vec U_{\ip{ns}+1,m},\ldots,\vec U_{\ip{nt},m}$ instead of $\vec U_{\ip{ns}+1},\ldots,\vec U_{\ip{nt}}$ in \eqref{Bn}.\\


The Proposition~\ref{prop1} can be seen as a corollary of the following Lemma:
\begin{lemma}\label{prop:Cnm}
Let $\vec X_1,\ldots,\vec X_n$ be $d$-dimensional independent random vectors with copula $C$, such that for a fixed integer $\bstar=\ip{nb}$, $b\in(0,1)$ known, the random vectors $\vec X_1,\ldots,\vec X_{\bstar}$ have m.c.d.f.s $F_1,\ldots,F_d$ and the random vectors $\vec X_{\bstar+1},\ldots,\vec X_n$ have m.c.d.f.s $F'_1,\ldots,F'_d$. Under Condition~\ref{Cj1}, 
$$
\sup_{(s,t,\vec u)\in\Delta\times[0,1]^d}|\C_{n,m}(s,t,\vec u) - \tilde\C_{n,m}(s,t,\vec u)|\p 0.
$$
consequently,
$$
\C_{n,m}\leadsto \C_C \qquad \text{ in } \ell^{\infty}(\Delta \times[0,1]^d).
$$
\end{lemma}
\newproof{Cnm}{Proof of Lemma \ref{prop:Cnm}}
\begin{Cnm}
We will demonstrated the 3 following convergences:
\begin{align}
\label{SUP5}&\sup_{(s,t,\vec u)\in(\Delta\cap[0,b]^2)\times[0,1]^d}|\C_{n,m}(s,t,\vec u) - \tilde\C_{n,m}(s,t,\vec u)|\p 0,\\
\label{SUP6}&\sup_{(s,t,\vec u)\in(\Delta\cap[b,1]^2)\times[0,1]^d}|\C_{n,m}(s,t,\vec u) - \tilde\C_{n,m}(s,t,\vec u)|\p 0,\\
\label{SUP7}&\sup_{(s,t,\vec u)\in[0,b]\times[b,1]\times[0,1]^d}|\C_{n,m}(s,t,\vec u) - \tilde\C_{n,m}(s,t,\vec u)|\p 0.
\end{align}
\end{Cnm}
On $(\Delta\cap[0,b]^2)\times[0,1]^d$ or $(\Delta\cap[b,1]^2)\times[0,1]^d$, we have $\C_{n,m}(s,t,\vec u)=\C_{n}(s,t,\vec u)$ and $\tilde \C_{n,m}(s,t,\vec u)=\tilde \C_{n}(s,t,\vec u)$,
hence the desirated convergence follow from \eqref{buckojseg}:
\begin{align}\label{cnmstar}
&\sup_{(s,t,\vec u)\in(\Delta\cap[0,b]^2)\times[0,1]^d}|\C_{n,m}(s,t,\vec u) - \tilde\C_{n,m}(s,t,\vec u)|\\
\nonumber &=\sup_{(s,t,\vec u)\in(\Delta\cap[0,b]^2)\times[0,1]^d}|\C_{n}(s,t,\vec u) - \tilde\C_{n}(s,t,\vec u)|\\
\nonumber &=\sup_{(s,t,\vec u)\in(\Delta\cap[0,b]^2)\times[0,1]^d}|\C^*_{n}(s,t,\vec u) - \tilde\C^*_{n}(s,t,\vec u)|\\
\nonumber &\leq\sup_{(s,t,\vec u)\in\Delta\times[0,1]^d}|\C^*_{n}(s,t,\vec u) - \tilde\C^*_{n}(s,t,\vec u)|\p0,
\end{align}
where $\C^*_{n}$ and $\tilde\C^*_{n}$ are the processes $\C_{n}$ and $\tilde\C_{n}$ constructed from $\vec X_1,\ldots,\vec X_{\ip{nb}},\vec X^*_{\ip{nb}+1},\ldots,\vec X^*_n$ with $\vec X^*_{\ip{nb}+1},\ldots,\vec X^*_n$ such as $\vec X_{1},\ldots,\vec X^*_{n}$ are i.i.d.

Because by construction, for all $(s,t,\vec u)\in[0,b]\times[b,1]\times[0,1]^d$ we have $\C_{n,m}(s,t,\vec u) = \C_{n}(s,b,\vec u)+\C_{n}(b,t,\vec u)$ and $\tilde \C_{n,m}(s,t,\vec u) = \tilde\C_{n}(s,b,\vec u)+\tilde\C_{n}(b,t,\vec u)$, the suprema in \eqref{SUP7} is bounded by
$$
\sup_{(s,\vec u)\in[0,b]\times[0,1]^d}|\C_{n}(s,b,\vec u) - \tilde\C_{n}(s,b,\vec u)| + \sup_{(t,\vec u)\in[b,1]\times[0,1]^d}|\C_{n}(b,t,\vec u) - \tilde\C_{n}(b,t,\vec u)|
$$

using the triangle inequality, hence the same argumentation as \eqref{cnmstar} allows us to conclude.\\

Using the fact that for any $m\in\{1,\ldots,n-1\}$, the random vectors $\vec U_{1,m},\ldots,\vec U_{n,m}$ are i.i.d., we obtain the weak limit of $(s,t,\vec u)\mapsto \B_{n,m}(s,t,\vec u)$ in $\ell^{\infty}(\Delta\times[0,1]^d)$ to $(s,t,\vec u)\mapsto\Z_C(t,\vec u) - \Z_C(s,\vec u)$  \citep[see for example the Theorem 2.12.1 of][]{vanWel96} hence the process $(s,\vec u)\mapsto\D_{n,\bstar}(s,\vec u)$ converges weakly in $\ell^{\infty}([0,1]^{d+1})$ to $\C^0_C(s,\vec u) - s\C^0_C(1,\vec u)$.
\end{pop}

\newproof{Cnmboot}{Proof of Proposition \ref{prop:Cnmboot}}
\begin{Cnmboot}
Proceeding as the proof of Proposition 4.3 in \cite{BucKojRohSeg14}, for $\beta\in\{1,\ldots,B\}$ and $(s,t,\vec u)\in\Delta\times[0,1]^d$, put
\begin{align*}
&\B^{(\beta)}_{n,m}(s,t,\vec u) = \frac{1}{\sqrt{n}}\sum_{i=\ip{ns}+1}^{\ip{nt}}\xi_i^{(\beta)}\{\1(\vec U_{i,m}\leq \vec u) - C( \vec u)\},\\
&\C^{(\beta)}_{n,m}(s,t,\vec u) = \B^{(\beta)}_{n,m}(s,t,\vec u)-\sum_{j=1}^d\dot C_{j}(\vec u)\B^{(\beta)}_{n,m}(s,t,\vec u^{\{j\}}).
\end{align*}
and the similar versions $\B^{(\beta)}_{n}$ and $\C^{(\beta)}_{n}$ based on $\vec U_{1},\ldots,\vec U_{n}$ instead of $\vec U_{1,m},\ldots,\vec U_{n,m}$. 
From Lemma~\ref{prop:Cnm}, the fact that $\vec U_{1,m},\ldots,\vec U_{n,m}$ are i.i.d. and using the Theorem 2.1 in \cite{BucKoj14} and the continuous mapping theorem, we have that
$$
\left(\C_{n,m},\C^{(1)}_{n,m},\ldots,\C^{(B)}_{n,m}\right)\leadsto \left(\C_C, \C^{(1)}_C,\ldots,\C^{(B)}_C\right)
$$
in $\ell^{\infty}(\Delta\times[0,1]^d)^{B+1}$. With this result, it is sufficient to demonstrate the following convergences: for $\beta=1,\ldots,B$,
\begin{align*}
&\sup_{(s,\vec u)\in[0,b]\times[0,1]^d}|\check \C^{(\beta)}_{n,m}(s,b,\vec u) - \C^{(\beta)}_{n,m}(s,b,\vec u)|\p 0,\\
&\sup_{(s,\vec u)\in[b,1]\times[0,1]^d}|\check \C^{(\beta)}_{n,m}(b,s,\vec u) - \C^{(\beta)}_{n,m}(b,s,\vec u)|\p 0,\\
\sup_{(s,\vec u)\in[0,b]\times[0,1]^d}|\check \C^{(\beta)}_{n,m}(0,s,\vec u) - 
&\C^{(\beta)}_{n,m}(0,s,\vec u)|\p 0,\qquad \sup_{(s,\vec u)\in[b,1]\times[0,1]^d}|\check \C^{(\beta)}_{n,m}(s,1,\vec u) - \C^{(\beta)}_{n,m}(s,1,\vec u)|\p 0,\\
\sup_{\vec u\in[0,1]^d}|\check \C^{(\beta)}_{n,m}(0,b,\vec u) - 
&\C^{(\beta)}_{n,m}(0,b,\vec u)|\p 0,\qquad \sup_{\vec u\in[0,1]^d}|\check \C^{(\beta)}_{n,m}(b,1,\vec u) - \C^{(\beta)}_{n,m}(b,1,\vec u)|\p 0.
\end{align*}

In fact, using the same argumentation as \eqref{cnmstar} and term (B.2) appearing in \cite{BucKojRohSeg14}, the previous convergences are automatically verified.

\end{Cnmboot}

\section{Simulation study}
\end{appendix}

\end{document}